\documentclass[12pt]{article}
\usepackage{geometry}
\usepackage{graphicx}
\usepackage{amssymb}
\usepackage{amsmath}
\usepackage{amsthm}
\usepackage{color}

\newtheorem{lemma}{Lemma}
\newtheorem{theorem}{Theorem}
\newtheorem{cor}{Corollary}

\newcommand{\R}{\mathbb R}
\numberwithin{equation}{section}

\title{Evaluating quasilocal energy and solving optimal embedding equation at null infinity}
\author{PoNing Chen, Mu-Tao Wang, and Shing-Tung Yau}

\begin{document}

\date{February 3, 2010, revised September 29, 2010}
\maketitle
\begin{abstract}
\footnote{M.-T. Wang is supported by NSF grant DMS 0904281 and S.-T. Yau is supported by NSF
grant PHY-0714648.} We study the limit of quasilocal energy defined in \cite{wy1} and \cite{wy2} for a family of spacelike 2-surfaces approaching null infinity of an asymptotically flat spacetime.
It is shown that Lorentzian symmetry is recovered and an energy-momentum 4-vector is obtained.
 In particular, the result is consistent with the Bondi--Sachs energy-momentum at a retarded time. The quasilocal mass in \cite{wy1} and \cite{wy2} is defined by minimizing quasilocal energy among admissible isometric embeddings and observers. The solvability of the Euler-Lagrange equation for this variational problem is also discussed in both the asymptotically flat and asymptotically  null cases. Assuming analyticity, the equation can be solved and the solution is locally minimizing in all orders. In particular, this produces an optimal reference hypersurface in the Minkowski space for the spatial or null exterior region of an asymptotically flat spacetime.
\end{abstract}

\section{Introduction}

This is a continuation of \cite{wy3} in which the spatial limit of the new quasilocal energy defined in \cite{wy1} and \cite{wy2} is analyzed. In the present article, we address the question of the null limit in Bondi--Sachs coordinates for an asymptotically flat spacetime.
Let $N$ be a spacetime with metric $g_{\alpha \beta}$ in Bondi--Sachs coordinates given by
$$
-UVdw^2-2Udwdr+\sigma_{ab}(dx^a+W^adw)(dx^b+W^bdw)\,\,\, a, b=2,3
$$
where
$$ W^a =O(r^{-2}),$$
$$U=1-\frac{X^2+Y^2}{2r^2}+o(r^{-2}),$$
$$V=1-\frac{2m}{r}+o(r^{-1})$$
and the metric $\sigma_{ab}$ is given by
\[ \left( \begin{array}{ccc}
r^2+2Xr+2(X^2+Y^2) & -2Yr \sin \theta  \\
 -2Yr \sin \theta & \sin^2 \theta[r^2-2Xr+2(X^2+Y^2)]
\end{array} \right)\] with
\[\det \sigma_{ab}=r^4\sin^2\theta.\]

The inverse of the metric $g_{\alpha\beta}$ is
\[g^{ww}=g^{wa}=0, g^{wr}=-U^{-1}, g^{rr}=U^{-1} V, g^{ra}=U^{-1} W^a, \text{\,and\,} g^{ab}=\sigma^{ab}.\]

Throughout the paper, coordinates are labeled by $x^0=w, x^1=r, x^2=\theta, x^3=\phi$ and the indexes are for $\alpha, \beta, \gamma \cdots =0, 1,2 ,3$,
$i, j, k \cdots=1, 2, 3$, and  $a, b \cdots=2,3$.

At a retarded time $w=c$, the Bondi--Sachs energy-momentum vector (\cite{bvm} \cite{s}) is defined as
\begin{equation}
(E, P_1, P_2, P_3)=\frac{1}{8 \pi} (\int_{S^2} 2m dS^2,   \int_{S^2} 2m\tilde{X}_1dS^2, \int_{S^2} 2m\tilde{X}_2dS^2, \int_{S^2} 2m\tilde{X}_3dS^2)
\end{equation} where $m=m(c, \theta, \phi)$ is the mass aspect function in the expansion of $V$, and
$\tilde{X}_i$, $i=1,2,3$ are the three eigenfunctions $\sin\theta \sin\phi$, $\sin\theta \cos\phi$ and $\cos\theta $ of the Laplace operator $\widetilde{\Delta}$ on $S^2$ with eigenvalue $-2$.

We recall that given a spacelike 2-surface $\Sigma$ in a spacetime, a quasilocal energy $E(\Sigma, X, T_0$) is defined in \cite{wy1}, \cite{wy2} with respect to an isometric embedding $X:\Sigma\rightarrow \R^{3,1}$ and a constant future timelike vector $T_0\in \R^{3,1}$.  For a family of surfaces $\Sigma_r$ and a family of isometric embeddings $X_r$ of $\Sigma_r$ into $\R^{3,1}$, the limit  of $E(\Sigma_r, X_r, T_0)$ is evaluated in \cite[Theorem 2.1]{wy3} under the assumption that
\begin{equation} \label{ratio}
\lim_{r \to \infty} \frac{|H_0|}{|H|}=1
\end{equation}
where $H$ and $ H_0$ are spacelike mean curvature vectors of $\Sigma_r$ in $N$ and  the image of $X_r$ in $\R^{3,1}$, respectively.
 In fact, the limit of $E(\Sigma_r, X_r, T_0)$ with respect to a constant future timelike vector $T_0 \in \R^{3,1}$ is given by
\begin{equation}\label{linear}
\lim_{r\to \infty}\frac{1}{8\pi}  \int_{\Sigma_{r}} \left[-\langle T_0,  \frac{{J}_0}{|H_0|}\rangle (|H_0|-|H|) - \langle\nabla^{\R^{3,1}}_{\nabla\tau} \frac{{J}_0}{|H_0|}, \frac{H_0}{|H_0|}\rangle+ \langle\nabla^{N}_{\nabla\tau} \frac{{J}}{|H|}, \frac{H}{|H|}\rangle \right] d \Sigma_r
\end{equation} where $\tau=-\langle T_0, X_r\rangle$ is the time function with respect to $T_0$, and $J_0$  and $J$ are the future timelike normal vectors dual to $H_0$ and $H$. This expression is linear in $T_0$ and defines an energy-momentum 4-vector at infinity.

In this article, we consider a family of 2-surfaces $\Sigma_r$ on a null cone $w=c$ as $r$ goes to infinity in Bondi--Sachs coordinates.
The limit of the quasilocal energy is first computed with respect to isometric embeddings $X_r$ into $\R^3$ which are essentially unique and satisfy (\ref{ratio}). We show in particular, \begin{equation}\lim_{r\to \infty}\frac{1}{8\pi}  \int_{\Sigma_{r}} (|H_0|-|H|) d\Sigma_r=E,\text{ and } \lim_{r\to \infty}\frac{1}{8\pi} \int_{\Sigma_{r}}\langle\nabla^{N}_{-\nabla X_i} \frac{{J}}{|H|}, \frac{H}{|H|}\rangle d \Sigma_r=P_i\end{equation} where $(X_1, X_2, X_3)$ are the coordinate functions of the isometric embedding $X_r$ into $\R^3$. We remark that exactly the same limit expression on coordinate spheres of asymptotically flat hypersurface gives the ADM energy-momentum in \cite{wy3}. The computation is stable with respect to any $O(1)$ perturbation of $X_r$ in $\R^{3,1}$ and is equivariant with respect to Lorentzian transformation acting on $X_r$.

In \cite{wy1} and \cite{wy2}, the quasilocal mass of a 2-surface $\Sigma$ is defined to be the minimum of $E(\Sigma, X, T_0)$ among all admissible pairs $(X, T_0)$ and the Euler-Lagrange equation is derived for an optimal isometric embedding. In the last section, we show that an analytic solution of the optimal isometric embedding equation can be obtained as an $O(1)$ perturbation of embeddings into a boosted totally geodesic slice in $\R^{3,1}$ whose timelike normal is in the direction of the total energy-momentum 4-vector. This solution locally minimizes the quasilocal energy.

 Brown--Lau--York \cite{blk} and Lau \cite{la} compute the null limit of the Brown--York energy and we compare our calculation with theirs in the following:

1) Brown--York mass is gauge dependent. After fixing a reference isometric embedding (either to flat $ \R^3$ \cite{blk} or to the null cone in $\R^{3,1}$\cite{la} ), a gauge is chosen arbitrarily so that the limit of the mass coincide with the Bondi mass.
In contrast, in our case, once a reference isometric embedding is picked, the quasilocal energy is determined by the canonical gauge condition (Eq (1.1) in \cite{wy3}). Our calculation is robust with respect to the choice of reference isometric embedding. In particular, the reference family can be arranged to be asymptotically flat or asymptotically null in $\R^{3,1}$.

2) In \cite{blk}, the momentum part came from the smear energy while in our case, the momentum part came from the connection one-form associated with the mean curvature gauge. This one form gives the right momentum contribution in the asymptotically flat case as well (see \cite{wy3}).

3) In \cite{blk}, the energy and momentum are defined separately. In our case, the Lorentzian symmetric is recovered at infinity and the energy-momentum form a covariant 4-(co)vector. We show that this (co)vector is equivariant with respect to the reference isometric embeddings into $\R^{3,1}$.

Acknowledgement: Part of the work is done while the authors are visiting the Taida Institute for Mathematical Sciences in Taipei, Taiwan.
%CH2
%CH2
%CH2
\section{The geometry of 2-surface $\Sigma_r$ in Bondi coordinates}
Let $N$ be an asymptotically flat spacetime with Bondi--Sachs coordinates.
Let $\Sigma_r$ be the 2-surface defined by $w=c$ and a fixed $r$. In this section, we compute the mean curvature vector $ H$ of $\Sigma_r$ in $N$ and the connection one-form of the normal bundle of $\Sigma_r$ in the mean curvature gauge. Denote $W_a=\sigma_{ab} W^b$ and let $\delta^a W_a$ be the divergence of the 1-form $W_a$ on $\Sigma_r$ with respect to the induced metric $\sigma_{ab}$.
\begin{lemma}
Let $\Sigma_r$ be the 2-surface defined by $w=c$ and a fixed $r$. The mean curvature vector $H$ of $\Sigma_r$ in $N$ is given by
\begin{equation} \label{eq:H1}
 H= \frac{1}{U} \left[ \frac{2}{r} ({\frac{\partial}{\partial w}} -W^{a}{\frac{\partial}{\partial x^a}} ) - (\frac{2V}{r}+ \delta^a W_a ){\frac{\partial}{\partial r}}\right].
\end{equation}
In particular, $H$ is spacelike when $r$ is large enough with
\begin{equation}\label{eq:H2}
|H|^2 = \frac{4}{Ur} (\frac{V}{r} +  \delta^a W_a).
\end{equation}  Suppose $J$ is the future timelike normal vector dual to $H$, then
\begin{equation}\label{eq:H3}  \langle \nabla _{\frac{\partial}{\partial x^b}} {J}, {H} \rangle =
                            \frac{2}{rU} \partial_{b} (\frac{V}{r}+\delta ^a W_a)-\frac{2}{r U^2}(\frac{V}{r}+{\delta^a W_a})\sigma_{cb}\partial_r W^c.\end{equation}

\end{lemma}
\begin{proof}
By definition, we have
\begin{align*}
 H= {} & \sigma^{ab} (\nabla _{\frac{\partial}{\partial x^a}}{\frac{\partial}{\partial x^b}} -  (\nabla _{\frac{\partial}{\partial x^a}}{\frac{\partial}{\partial x^b}}) ^{T})\\
   = {} & \sigma^{ab} (\Gamma_{ab}^r {\frac{\partial}{\partial r}} + \Gamma_{ab}^w {\frac{\partial}{\partial w}} + (\Gamma_{ab}^c -\langle \nabla _{\frac{\partial}{\partial x^a}}{\frac{\partial}{\partial x^b}},{\frac{\partial}{\partial x^d}}\rangle  \sigma^{dc} ){\frac{\partial}{\partial x^c}}  ).
 \end{align*}
The last coefficient can be computed explicitly as
 \begin{align*}
   &  \Gamma_{ab}^c -\langle \nabla _{\frac{\partial}{\partial x^a}}{\frac{\partial}{\partial x^b}},{\frac{\partial}{\partial x^d}}\rangle  \sigma^{dc}  \\
 ={} &  \Gamma_{ab}^c -(\Gamma_{ab}^r g_{rd}+\Gamma^{w}_{ab} g_{wd}+\Gamma_{ab}^e\sigma_{ed}) \sigma^{dc}\\
  ={} & - \Gamma_{ab}^r g_{rd}\sigma^{dc} - \Gamma_{ab} ^w g_{wd} g^{dc} \\
  ={} & \Gamma_{ab} ^w g_{wr} g^{rc}.
 \end{align*}
 Thus $$
 H= \sigma^{ab} (\Gamma_{ab}^r {\frac{\partial}{\partial r}} + \Gamma_{ab}^w {\frac{\partial}{\partial w}} + \Gamma_{ab} ^w g_{wr} g^{rc} {\frac{\partial}{\partial x^c}} ).
$$
The  relevant Christoffel symbols of $g_{\alpha \beta}$ are given by
$$
\Gamma^w_{ab} = \frac{1}{2} U^{-1} \partial_r \sigma _{ab},
$$ and
\[
\Gamma^r_{ab}=- \frac{1}{2} U^{-1} [ \partial_b W_a + \partial_a W_b -  \partial_w \sigma _{ab} +V \partial_r \sigma _{ab}- 2 \gamma^d_{ab}W_d] \] where $\gamma_{ab}^d$
is the Christoffel symbol of the metric $\sigma_{ab}$.
When tracing with $\sigma^{ab}$, we notice that
\begin{equation}\label{eq:diff}
\sigma^{ab}  \partial_{{\alpha}} \sigma _{ab} = \partial_{{\alpha}} (\ln \det \sigma_{ab}) = \partial_{{\alpha}} \ln (r^4 \sin^2 \theta) .
\end{equation}
Thus, we obtain equations (\ref{eq:H1}) and (\ref{eq:H2}).
To compute the connection one-form, we rewrite equation (\ref{eq:H1}) as
\begin{align*}
U {H} ={} & -(\frac{2V}{r}+{\delta ^a W_a}) \frac{\partial}{\partial r} + \frac{2}{r}(  \frac{\partial}{\partial w} -W^c \frac{\partial}{\partial x^c}) \\
       = {} &  -(\frac{V}{r}+{\delta ^a W_a}) \frac{\partial}{\partial r} + \frac{2}{r}(  \frac{\partial}{\partial w} -W^c \frac{\partial}{\partial x^c} -\frac{V}{2} \frac{\partial}{\partial r} )
\end{align*}
where $ \frac{\partial}{\partial r}$ and $  \frac{\partial}{\partial w} -W^c \frac{\partial}{\partial x^c} -\frac{V}{2} \frac{\partial}{\partial r}$ are null vectors.

Thus we have
$$
U {J} =  -(\frac{V}{r}+{\delta ^a W_a}) \frac{\partial}{\partial r} - \frac{2}{r}(  \frac{\partial}{\partial w} -W^c \frac{\partial}{\partial x^c} -\frac{V}{2} \frac{\partial}{\partial r} ).
$$
For simplicity, let's denote  $ \frac{\partial}{\partial r}$ and $  \frac{\partial}{\partial w} -W^c \frac{\partial}{\partial x^c} -\frac{V}{2} \frac{\partial}{\partial r}$  by $ \vec n_1$ and $\vec n_2$ and the coefficients $\frac{V}{r}+\delta ^a W_a$ and $\frac{2}{r}$ by $x$ and $y$ in the following computation. Then,
\begin{equation}\label{JH}\begin{split}
 \langle \nabla _{\frac{\partial}{\partial x^b}} J, H \rangle ={} & U^{-2}  \langle\nabla _{\frac{\partial}{\partial x^b}} x\vec{n}_1+y\vec{n}_2, x
 \vec{n}_1-y\vec{n}_2\rangle\\
              ={} &U^{-2} \{ [ (\partial_{b}x)(-y) + (\partial_{b}y)(x)] \langle \vec{n}_1, \vec{n}_2 \rangle -xy( \langle \nabla _{\frac{\partial}{\partial x^b}}
               \vec{n}_1, \vec{n}_2
             \rangle- \langle \nabla _{\frac{\partial}{\partial x^b}} \vec{n}_2, \vec{n}_1 \rangle) \}\\
               ={} &U^{-2} [ (\partial_{b}x)(-y) \langle \vec{n}_1, \vec{n}_2\rangle  -2xy\langle \nabla _{\frac{\partial}{\partial x^b}} \vec{n}_1, \vec{n}_2
             \rangle +xy  \partial_b \langle \vec{n}_2, \vec{n}_1 \rangle ]
 \end{split}\end{equation}
On the other hand,
\[ \langle \nabla _{\frac{\partial}{\partial x^b}} \vec{n}_1, \vec{n}_2 \rangle=\langle   \Gamma^r_{br}
                \frac{\partial}{\partial r} ,  \frac{\partial}{\partial w} -W^d \frac{\partial}{\partial x^d}
                 -\frac{V}{2} \frac{\partial}{\partial r} \rangle= -U\Gamma_{br}^r\] because $\Gamma_{br}^w=0$ by direct computation, and $\frac{\partial}{\partial x^c}$ is perpendicular to the null normal $  \frac{\partial}{\partial w} -W^d \frac{\partial}{\partial x^d}-\frac{V}{2} \frac{\partial}{\partial r} $.
Substitute in $\Gamma_{br}^r=\frac{1}{2} U^{-1}\partial_b U-\frac{1}{2}U^{-1} \sigma_{bc} \partial_r W^c$, and we obtain
\begin{equation} \label{n1n2}\langle \nabla _{\frac{\partial}{\partial x^b}} \vec{n}_1, \vec{n}_2 \rangle=-\frac{1}{2}\partial_b U+\frac{1}{2}\sigma_{bc}\partial_r W^c.\end{equation}

Plug \eqref{n1n2} and $\langle \vec{n}_1,\vec{n}_2\rangle=-U$ into \eqref{JH}, we derive
\begin{align*}  \langle \nabla _{\frac{\partial}{\partial x^b}} {J}, {H} \rangle =&
                              U^{-2}y[(\partial_b x) U-x\sigma_{bc}\partial_r W_c],\end{align*}
and \eqref{eq:H3} follows in view of the definitions of $x$ and $y$.
\end{proof}
%CH3
%CH3
%CH3
\section{Limit of quasilocal energy}
In this section, we compute the limit of quasilocal energy with respect to a family of isometric embeddings $X_r$ of $\Sigma_r$ as an $O(1)$ perturbations of a boosted totally geodesic slice in $\R^{3,1}$. First we quote the following lemma whose proof can be found in \cite{fst}:

\begin{lemma}\label{isom_R^3}
Let $\sigma^r_{ab}$ be a family of metrics on $\Sigma_r\simeq  S^2$ with $\sigma^r_{ab}=r^2\tilde{\sigma}_{ab}+O(r)$ in which $\tilde{\sigma}_{ab}$ is the standard round metric on $S^2$. Let $X_r=(X_1, X_2, X_3)$ be the isometric embedding into $\R^3$ for $r$ large and $H_0$ be the mean curvature of $X_r$. Then
\[|{H}_0| = \frac{2}{r} +O(r^{-2}) \text{  and  }\int_{\Sigma_r} |{H}_0| d\Sigma_r =  4 \pi r + \frac{Area(\Sigma_r)}{r}+O(r^{-1}).
\]
\end{lemma}
 We note that up to an isometry of $\R^3$, $X_r$ can be arranged so that the coordinate functions satisfy $X_i=r\tilde{X}_i+ O(1)$.
\begin{theorem}\label{maintheorem}
Let $\Sigma_r$ be the 2-surface defined by $w=c$ and a fixed $r$ in an asymptotically flat spacetime with Bondi--Sachs coordinates. Suppose $X_r$ is the (unique) family of isometric embeddings of $\Sigma_r$ into $\R^{3}$ for $r$ large, the limit of quasilocal energy with respect to $T_0=(\sqrt{1+|a|^2}, a^1,a^2,a^3 )$ is
\begin{equation}
\lim_{r\rightarrow\infty} E(\Sigma_r, X_r, T_0)=\frac{1}{8 \pi} \int_{S^2} 2m (\sqrt{1+|a|^2}+ a^i\tilde{X}_i)dS^2.
\end{equation}
\end{theorem}
\begin{proof}
Let $(0,X_1,X_2,X_3)$ be the isometric embedding $X_r$ of $\Sigma_r$ into $\R^3\subset \R^{3,1}$. In this case, $\frac{J_0}{|H_0|}$ is simply the vector $(1,0,0,0)$. By the assumption on $\sigma_{ab}$ we can apply Lemma \ref{isom_R^3} and
 $$
  \int_{\Sigma_r} |H_0| d\Sigma_r =  8 \pi r +O(r^{-1}).
 $$ On the other hand, from equation (\ref{eq:H2}) and the expansion for $V$, we obtain
\begin{equation} \label{normH}
|H| = \frac{2}{r} -\frac{2m}{r^2} + \delta^a W_a +O(r^{-3})
\end{equation} and thus
 $$
  \int_{\Sigma_r} |H| d\Sigma_r =   8 \pi r  - \int_{S^2}  2m \, dS^2 + O(r^{-1}).
 $$
 Next we compute the physical hamiltonian
 \[ \frac{1}{8 \pi} \int_{\Sigma_r}  \langle \nabla^N_{\nabla \tau} \frac{J}{|H|} , \frac{ H}{|H|}   \rangle  d\Sigma_r=-a^i \frac{1}{8 \pi} \int_{\Sigma_r}  \langle \nabla^N_{\nabla X_i} \frac{J}{|H|} , \frac{ H}{|H|}   \rangle  d\Sigma_r.\]
From equation (\ref{eq:H3}) and the asymptotic expansions of $V$ and $W^a$ , we derive
\[
\langle \nabla _{\frac{\partial}{\partial x^b}} {J}, {H} \rangle =\frac{2}{r}[\partial_{ b} (\delta ^a W_a-\frac{2m}{r^2})] + \frac{4}{r^3} W_b+O(r^{-4}).\]

Let $V$ denote the connection one-from $\langle \nabla^N \frac{J}{|H|} , \frac{ H}{|H|}   \rangle$. From the above computation,
\[  div_{\Sigma_r} V =  \frac{1}{r}[\frac{1}{2}(\widetilde{\Delta}+2)({\delta^a}{W_a} -\frac{2m}{r^2}) +  \frac{2m}{r^2}] +O(r^{-4}). \]

The limit of $\int_{\Sigma_r}  \langle \nabla^N_{\nabla X_i} \frac{J}{|H|} , \frac{ H}{|H|}   \rangle  d\Sigma_r$ as $r\rightarrow \infty$ is thus the same as
\[ \lim_{r\rightarrow \infty}\int_{\Sigma_r} X_i div_{\Sigma_r} V d\Sigma_r  =  \int_{S^2} \tilde X_i  [\frac{1}{2}(\widetilde{\Delta}+2)(\widetilde{\delta^a}\widetilde{W_a} -2m)+  2m ]dS^2 =  \int_{S^2} \tilde X_i 2m \, dS^2\]

In this case, the reference Hamiltonian term is zero as $\frac{J_0}{|H_0|}$ is a constant vector. In view of expression (\ref{linear}), the theorem is proved.

\end{proof}

Next we show that the limit of the quasilocal energy is invariant under any $O(1)$ perturbations of embeddings into totally geodesic $\R^3$ and that it is Lorentzian equivariant.
\begin{cor}
Suppose $X_r=(\tau_0, X_1, X_2, X_3)$ is a family of isometric embeddings of $\Sigma_r$ into $\R^{3,1}$ with $ \tau_0 =  \tilde \tau_0+O(r^{-1})$ for some function $\tilde \tau_0$ on $S^2$.
Then we still have
\[\lim_{r\rightarrow\infty} E(\Sigma_r, X_r, T_0)=\frac{1}{8 \pi} \int_{S^2} 2m (\sqrt{1+|a|^2}+ a^i\tilde{X}_i)dS^2.\]
\end{cor}

\begin{proof}
Let $\hat{X}_r$ be the embedding of $\Sigma_r$ by projecting $X_r$ onto $\R^{3}$ which is given by, $(0, X_1, X_2, X_3)$. It is not hard to check that the induced metric by the embedding $\hat{X}_r$ agrees with the standard round metric of radius $r$ up to the top order term and its area  agrees with that of the standard round metric of radius $r$ up to the second order term.
The mean curvature  of the  embedding $\hat{X}_r$ is then $ (0, \widehat\Delta X_1, \widehat\Delta X_2 , \widehat\Delta X_3)$. By Lemma \ref{isom_R^3}, the mean curvature $\hat{H}_0$ satisfies
 \begin{equation}\label{hat_H}|\hat{H}_0| = \frac{2}{r} +O(r^{-2}), \text{  and  }\int_{\Sigma_r} |\hat{H}_0| d\Sigma_r =  8 \pi r + O(r^{-1}).
\end{equation}
The mean curvature $ H_0$ of $X_r$ is given by
\[ (\Delta \tau_0, \Delta X_1, \Delta X_2 , \Delta X_3).    \]
The difference between $\sigma_{ab}$ and $\hat{\sigma}_{ab}$ is of order
\[  \sigma_{ab} -\hat{\sigma}_{ab} =O(1) \text{ and } \sigma^{ab} -\hat{\sigma}^{ab} =O(r^{-4}). \]
As a result, the difference between the two laplace operators is of order
\[   \widehat\Delta X_i - \Delta X_i =O(r^{-3}).  \]
Hence $|\hat{H}_0|^2 -|H_0|^2 =O(r^{-4})$, and thus $|\hat{H}_0| -|H_0| =O(r^{-3})$.

By equation (\ref{linear}), the  limit of quasilocal energy with respect to the embedding $(\tau_0, X_1, X_2, X_3)$ is thus
\begin{align*}
  \frac{1}{8 \pi}\lim_{ r\to \infty}  \int_{\Sigma_r} (|H_0| -|H|)  d\Sigma_r=\frac{1}{8 \pi} \lim_{ r\to \infty}  \int_{\Sigma_r} (|\hat{H}_0| -|H|) d\Sigma_r.
    \end{align*}
Unlike the previous case, $\frac{{J}_0}{|H_0|} $ is no longer a constant vector for such an isometric embedding $X_r$. However, the asymptotic expansion
$$
-\langle T_0,  \frac{{J}_0}{|H_0|} \rangle = \sqrt{1+|a|^2} +O(r^{-1})
$$ is valid and the energy component is the same as the  limit of quasilocal energy of the isometric embedding into $\R^3$ in view of (\ref{hat_H}).

 Next we compute the physical hamiltonian. Since the induced metric on the projection still agrees with the standard one up to lower order term, up to an isometry of $\R^3$, $X_i=r\tilde{X}_i+O(1)$.  The corresponding time function  $\tau$ is
$$  \tau =  - (\sum_i a_i X_i)  +  \tau_0\sqrt{1+\sum_i a_i^2} = -(\sum_i a_i \tilde X_i) r + O(1).$$
 Thus the physical hamiltonian remains the same.

 Lastly, we claim that the reference hamiltonian goes to $0$ as $r$ goes to infinity.
 \begin{lemma}\label{meancurvature}
Let $\sigma_r = r^2 \tilde \sigma +O(r)$ be a family of metrics on $S^2$.
Given an $O(1)$ time function $\tau_0$, let $ \Sigma^0_r $ be the images of the  isometric embedding $X_r$ into $\R^{3,1}$ determined by $\tau_0$.
Let $V_0$ be the vector dual to the one form $\langle  \nabla_{(\cdot) } \frac{ J_0 }{|H_0|}, \frac{H_0}{|H_0|} \rangle$ on $\Sigma^0_r$
then
\[ div_{\Sigma_r}V_0 =   \frac{1}{2r^3} \widetilde{\Delta}(\widetilde{\Delta}+2) \tau_0 +O(r^{-4}). \]
 \end{lemma}
 \begin{proof}
We need to compute $\frac{ J_0}{|H_0|}$ up to $O(r^{-2})$. For this purpose, it is enough to assume that
the embedding is $(\tau_0, r \tilde{X}_1, r \tilde{X}_2, r\tilde{X}_3 ) $.
Using $\sum_{i} \partial_c\tilde{X}_i\partial_b \tilde{X}_i=\tilde{\sigma}_{bc}$, we derive that a normal vector is
$$
(1, \frac{1}{r} \partial_a \tau_0 \partial_b \tilde X_i\tilde \sigma^{ab})
$$
where $\tilde \sigma$ is the standard metric on $S^2$. As $\sum_{i} \tilde{X}_i\partial_a \tilde{X}_i=0$, we check that $$
(1, \frac{1}{r} \partial_a \tau_0  \partial_b \tilde X_i\tilde \sigma^{ab}) +\frac{\Delta \tau_0}{|H_0|^2}  H_0
$$ is a normal vector perpendicular to $H_0$. Thus, $\frac{ J_0}{|H_0|}$, which is the unit normal perpendicular to $ H_0$ is, up to lower order, given by the same expression. As a result, we compute
$$
\langle  \nabla_{\frac{\partial}{\partial x^a} } \frac{ J_0 }{|H_0|}, H_0   \rangle=
\sum_i  \partial_a \left[\frac{1}{r}  (\partial_b \tau_0 \partial_c \tilde X_i)\tilde \sigma^{bc}+\frac{\Delta \tau_0}{|H_0|^2}\Delta X_i \right] \Delta X_i.
$$
The right hand side equals to
$$\left[\sum_i \frac{1}{r} (\partial_b \tau_0 )(\partial_a \partial_c \tilde X_i)\tilde \sigma^{bc}\Delta X_i\right] + \partial_a \Delta \tau_0 +O(r^{-3}).$$
Here one uses again that  $\sum_{i} \tilde{X}_i\partial_a \tilde{X}_i=0$ and thus from the first term, one has non-zero contribution only when the derivative $ \partial_a $ falls on $ \partial_c \tilde X_i$. For the second term, the leading term of $|H_0|$ is independent of $\theta$ and $\phi$ and $\sum _i (\Delta X_i) ^2 =|H_0|^2$ up to lower order terms. Thus one only has contribution when the derivative hits $\Delta \tau_0$. Direct computation using $\sum_i \tilde{X}_i \partial_a\partial_c\tilde{X}_i=-\tilde{\sigma}_{ac}$ shows that
$$
\sum_i \frac{1}{r}  (\partial_b \tau_0) (\partial_a\partial_c \tilde X_i)\tilde \sigma^{bc}\Delta X_i  = \frac{2}{r^2}\partial_a \tau_0.
$$
As a result,
\[ div_{\Sigma_r}V_0 =   \frac{1}{2r^3} \widetilde{\Delta}(\widetilde{\Delta}+2) \tau_0 +O(r^{-4}) \]
\end{proof}
Using the above lemma, the reference hamiltonian at infinity is
\[
\lim_{r\to\infty}\int_{\Sigma_r} X_i div_{\Sigma_r} V_0  d\Sigma=  \int_{S^2}  \frac{1}{2} \tilde X_i\widetilde{\Delta}(\widetilde{\Delta}+2) \tau_0 dS^2=0
\]
\end{proof}
\begin{cor}
 Suppose $X_r'$ is another family of isometric embeddings of $\sigma_r$ into $\R^{3,1}$ such that $X'_r=\tilde{L}_r X_r$ for some $X_r$ in the previous Corollary, and a family of Lorentzian transformation $\tilde{L}_r$ such that the limit of the $SO(3,1)$ part of $\tilde{L}_r$ converges to an $L_\infty$, then the energy-momentum 4-vector also transform by $L_\infty$.
\end{cor}

\begin{proof}
Both $|H_0|$ and the connection one form $\langle\nabla_{(\cdot)}^{\R^{3,1}} \frac{{J}_0}{|H_0|}, \frac{H_0}{|H_0|}\rangle$ are invariant under Lorentzian transformation, while $\langle T_0,  \frac{{J}_0}{|H_0|}\rangle$ and $\nabla\tau= - \nabla \langle T_0, X_r\rangle$ are Lorentzian equivariant.
\end{proof}

For example, if we take a family of isometric embedding $X_r$ into $\R^3$ and define $X'_r=X_r+r$, it is not hard to see that the hypersurface spanned by $X'_r$ is asymptotically null.
%CH4
%CH4
%CH4
\section{Optimal embedding equation}\label{opt}
The optimal embedding equation for minimizing the quasilocal energy is derived in \cite[Proposition 6.2]{wy2}. The equation reads
\begin{equation}\label{optimal}
-(\widehat{H}\hat{\sigma}^{ab} -\hat{\sigma}^{ac} \hat{\sigma}^{bd} \hat{h}_{cd})\frac{\nabla_b\nabla_a \tau}{\sqrt{1+|\nabla\tau|^2}}+ div_\Sigma (\frac{\nabla\tau}{\sqrt{1+|\nabla\tau|^2}} \cosh\theta|{H}|-\nabla\theta-V)=0\end{equation}
where $\sinh \theta =\frac{-\Delta \tau}{|{H}|\sqrt{1+|\nabla \tau|^2}} $ and $V$ is the connection one-form $\langle \nabla^N_{(\cdot)} \frac{J}{|H|}, \frac{H}{|H|}\rangle$. To solve for this equation, we start with data on the 2-surface $\Sigma$ given by $(\sigma_{ab}, |H|, V)$. Take a function $\tau$ on $\Sigma$ and consider the isometric embedding $\hat{X}:(\Sigma, \hat{\sigma}) \rightarrow \R^3$ with the metric $\hat{\sigma}_{ab}=\sigma_{ab}+\tau_a \tau_b$. $\widehat{H}$ and $\hat{h}_{ab}$ are the mean curvature and the second fundamental form of $\hat{\Sigma}$, the image of $\hat{X}$ in $\R^3$, respectively.

 Let $X:\Sigma\rightarrow \R^{3,1}$ be the embedding of the graph of $\tau$ over $\hat{\Sigma}$ in $\R^{3,1}$, and $\Sigma_0$ be the image of $X$ with the induced metric isometric to $\sigma$. The optimal isometric embedding equation can be written in terms of the geometry of $\Sigma_0$. In fact, the quasilocal energy of $\Sigma_0$ with respect to itself as a reference is zero and thus minimizing.
 Equation \eqref{optimal} is automatically true on $\Sigma_0$ and we deduce
\begin{equation}\label{optimal2}
-(\widehat{H}\hat{\sigma}^{ab} -\hat{\sigma}^{ac} \hat{\sigma}^{bd} \hat{h}_{cd})\frac{\nabla_b\nabla_a \tau}{\sqrt{1+|\nabla\tau|^2}}+ div_\Sigma (\frac{\nabla\tau}{\sqrt{1+|\nabla\tau|^2}} \cosh\theta_0|{H_0}|-\nabla\theta_0-V_0)=0
\end{equation} where $V_0$, $H_0$, $\theta_0$ are the corresponding data on $\Sigma_0$. Equation \eqref{optimal2} can be checked directly for spacelike surfaces in $\R^{3,1}$.

Subtracting equation (\ref{optimal}) from equation (\ref{optimal2}), equation (\ref{optimal}) is then equivalent to \begin{equation}\label{optimal2.5} div_{\Sigma} [ \frac{\nabla \tau }{\sqrt{1+ |\nabla \tau|^2}}(\cosh \theta |H| - \cosh \theta_0 |H_0|)- \nabla(\theta - \theta_0) -V  + V_0] =0.\end{equation}
By the definition of $\theta$ and $\theta_0$, we derive
\[ \cosh \theta |H| - \cosh \theta_0 |H_0| =\sqrt{|H|^2 +\frac{(\Delta \tau)^2}{1+ |\nabla \tau|^2}} - \sqrt{|H_0|^2 +\frac{(\Delta \tau)^2}{1+ |\nabla \tau|^2}}   \] and
\[  \sinh (\theta -\theta _0 ) =  \frac{\Delta \tau}{|H|| H_0|  \sqrt{1+ |\nabla \tau|^2}}(\sqrt{|H|^2 +\frac{(\Delta \tau)^2}{1+ |\nabla \tau|^2}} - \sqrt{|H_0|^2 +\frac{(\Delta \tau)^2}{1+ |\nabla \tau|^2}}). \]
Set

\begin{equation}\label{form_f} \begin{split}f &= \frac{\sqrt{|H_r|^2 +\frac{(\Delta \tau)^2}{1+ |\nabla \tau|^2}} - \sqrt{|H_0|^2 +\frac{(\Delta \tau)^2}{1+ |\nabla \tau|^2}} }{ \sqrt{1+ |\nabla \tau|^2}}\\
&=\frac{|H_r|^2-|H_0|^2}{\sqrt{|H_r|^2(1+|\nabla\tau|^2)+(\Delta\tau)^2}+\sqrt{|H_0|^2(1+|\nabla\tau|^2)+(\Delta\tau)^2}}. \end{split}\end{equation}and equation \eqref{optimal} is equivalent to
\begin{equation} \label{optimal3}
div_\Sigma(f \nabla \tau) - \Delta [ \sinh^{-1} (\frac{\Delta \tau f}{|H||H_0|})]-( div_\Sigma V -div_\Sigma V_0)=0.
\end{equation}

In this equation, $V$ are $|H|$ come from the physical data, and $V_0$ and $|H_0|$ only depend on the embedding $X:\Sigma\rightarrow \R^{3,1}$, while $\tau=-\langle X, T_0\rangle$ depends on both $X$ and $T_0$. Equation \eqref{optimal3} together with the isometric embedding equation $\langle dX, dX\rangle=\sigma$ form the optimal isometric embedding system.

We shall solve the system for a family of spacelike 2-surfaces at null or spatial infinity such that the family of isometric embeddings $X_r$ into $\R^{3,1}$ is of the form $X_r=B_r\hat{X}_r$ where $B_r$ is a family in $SO(3,1)$ and $\hat{X}_r$ is an $O(1)$ perturbation of isometric embeddings into $\R^3$.

We observe that momentum become an obstruction to solving the optimal embedding equation for $\hat{X}_r$ and then discuss how this can be resolved by boosting the embedding by $B_r$. The discussion covers the spatial infinity case discussed in \cite{wy3} as well. In the last subsection, we show the solution obtained is locally energy-minimizing up to lower order terms in $r$.
%CH4.1
%CH4.1
%CH4.1
\subsection{Embedding near $\R^3$}
In this subsection, we study the geometry of a family of isometric embeddings $\hat{X}_r$ that is near a totally geodesic $\R^3$ in $\R^{3,1}$.

\begin{lemma}\label{mean_curv} Suppose $\hat{X}_r$ is a family of isometric embeddings into $\R^{3,1}$ for a given family of metrics $\sigma_r=r^2\tilde{\sigma}+r\sigma^{(1)}+\sum_{k=0}^\infty r^{-k} \sigma^{(-k)}$ and
\begin{equation}\label{hat_X}\hat{X}_r= r\hat{X}^{(1)}+ \sum_{k=0}^\infty r^{-k} \hat{X}^{(-k)}\end{equation}
with $\hat{X}^{(1)}=\tilde{X}=(0, \tilde{X}_1, \tilde{X}_2, \tilde{X}_3)$, the standard embedding of $S^2$ into $\R^3$.
Denote by $\hat{\tau}^{(k)}$ the time component of $\hat{X}^{(k)}$ and by $G(\hat{\tau}^{(0)},\cdots, \hat{\tau}^{(-l)},\sigma) $
a term that depends on $\hat{\tau}^{(0)},\cdots, \hat{\tau}^{(-l)}$ and $\sigma$.  Then
\[|H_0|={2}r^{-1}+r^{-2} h_0^{(-2)}+\sum_{k=3}^\infty r^{-k} h_0^{(-k)}\]
where $h_0^{(-k)}=G(\hat{\tau}^{(0)},\cdots, \hat{\tau}^{(-k+3)},\sigma)$ for $k\geq 3$ and $h_0^{(-2)}$ depends only on $\sigma$.

\end{lemma}
\begin{proof} In the proof, we suppress the subscript $r$ and write $\hat{X}$ for $\hat{X}_r$ and $\sigma$ for $\sigma_r$.  From the expansion of $\sigma$, $\sigma=r^2\tilde{\sigma}+r\sigma^{(1)}+\cdots$, we deduce
\[2\sum_{i=1}^3 d \tilde{X}_i d\hat{X}_i^{(0)}=\sigma^{(1)}\]
This can be transformed into a linear elliptic equation for $\hat{X}_i^{(0)}$ which can be solved (see section 6 of \cite{n}).
In general,
\[2 \sum_{i=1}^3 d \tilde{X}_i d\hat{X}_i^{(-l-1)}+\sum_{m=0}^l \langle d\hat{X}^{(-m)}, d \hat{X}^{(-l+m)}\rangle=\sigma^{(-l)}\]
and thus $\hat{X}_i^{(-l-1)}$ is determined by $\hat{\tau}^{(0)},\cdots, \hat{\tau}^{(-l)}$ and $\sigma$.

Recall the mean curvature is given by $H_0={\Delta}\hat{X}$ where ${\Delta}$ is the Laplace operator with respect to $\sigma$.
For a function $g$ on $S^2$, we compute
\[\begin{split}\Delta g&=\sigma^{ab}(\partial_a\partial_b g-\gamma_{ab}^c\partial_c g)\\
&=(r^{-2}\tilde{\sigma}^{ab}+r^{-3} \sigma^{(-3) ab}+O(r^{-4}))(\widetilde{\nabla}^2_{ab}g-r^{-1} \gamma^{(-1)c}_{ab} \partial_c g+O(r^{-2}))\\
\end{split}\]
where $\gamma_{ab}^c=\tilde{\gamma}^{c}_{ab}+r^{-1}\gamma^{(-1)c}_{ab}+\cdots$ is the expansion of the  Christoffel symbol $\gamma_{ab}^c$ of
$\sigma_{ab}$ and $\widetilde{\nabla}^2_{ab}g=\partial_a\partial_b g-\tilde{\gamma}_{ab}^c\partial_c g$ is the Hessian of $g$ with respect to $\tilde{\sigma}_{ab}$.
Therefore, we obtain the following formula:
\begin{equation}\label{Delta_g}\Delta g=r^{-2}\widetilde{\Delta} g+r^{-3} (\sigma^{(-3) ab} \widetilde{\nabla}^2_{ab} g-\tilde{\sigma}^{ab}\gamma_{ab}^{(-1)c}\partial_c g)+O(r^{-4})\end{equation}
By \eqref{hat_X}, $H_0=\Delta \hat{X}$ has the following expansion:
\[ H_0=r\Delta \hat{X}^{(1)}+\sum_{k=0}^\infty r^{-k} \Delta\hat{X}^{(-k)}.\]

Since $\hat{X}^{(1)}=\tilde{X}$, the standard embedding of $S^2$, we compute
\[\Delta \tilde{X}=-2r^{-2}\tilde{X}+r^{-3}(-\tilde{\sigma}_{ab}\sigma^{(-3)ab}\tilde{X}-\tilde{\sigma}^{ab}{\gamma^{(-1)}}_{ab}^c\partial_c \tilde{X})+O(r^{-4})\]
where we use $\widetilde{\nabla}^2_{ab}\tilde{X}=-\tilde{\sigma}_{ab}\tilde{X}$.

Therefore, we obtain
\[H_0=\Delta \hat{X}= -2r^{-1}\tilde{X}+r^{-2}H_0^{(-2)}+\sum_{k=3}^\infty r^{-k}H_0^{(-k)}\]
where
\[H_0^{(-2)}=\widetilde{\Delta}\hat{X}^{(0)}-(\tilde{\sigma}_{ab}{\sigma^{(-3)}}^{ab}\tilde{X}+\tilde{\sigma}^{ab}{\gamma^{(-1)}}_{ab}^c\partial_c\tilde{X})\]
 and
\[H_0^{(-k)}=\widetilde{\Delta}\hat{X}^{(-k+2)}+G(\hat{\tau}^{(0)},\cdots \hat{\tau}^{(-k+3)},\sigma)\] for $k\geq 3$.

We compute
\[|H_0|^2=4r^{-2}+r^{-3}(2\widetilde{\Delta}\tilde{X}_i\widetilde{\Delta}\hat{X}_i^{(0)}+G(\sigma)) +\sum_{k=3}^\infty r^{-k-1}(2\sum_{i=1}^3\widetilde{\Delta}\tilde{X}_i\widetilde{\Delta} \hat{X}_i^{(-k+2)}+G(\hat{\tau}^{(0)},\cdots, \hat{\tau}^{(-k+3)},\sigma)  ).\] Since $\hat{X}_i^{(-k+2)}$ depends only on  $\hat{\tau}^{(0)},\cdots, \hat{\tau}^{(-k+3)}$ and $\hat{X}_i^{(0)}$ depends only on $\sigma$, we obtain
\[|H_0|^2={4}r^{-2}+r^{-3}G(\sigma)+\sum_{k=0}^\infty r^{-k-4} G(\hat{\tau}^{(0)},\cdots, \hat{\tau}^{(-k)},\sigma)\] and the expansion for $|H_0|$ follows by taking square root.
\end{proof}

\begin{lemma}\label{conn}
Under the same assumption as Lemma \ref{mean_curv}.
Let $V_0$ be the one form $\langle  \nabla_{(\cdot) } \frac{ J_0 }{|H_0|}, \frac{H_0}{|H_0|} \rangle$ on the image of $\hat{X}_r$
then
\[div_{\Sigma_r}V_0=r^{-3}[\frac{1}{2} \widetilde\Delta(\widetilde \Delta +2) \hat{\tau}^{(0)}]+  \sum_{k=4}^{\infty} r^{-k} v_0^{(-k)}  \]
where
\[v^{(-k)}=\frac{1}{2} \widetilde\Delta(\widetilde \Delta +2) \hat{\tau}^{(-k+3)} + G(\hat{\tau}^{(0)} , \cdots , \hat{\tau}^{(-k+4)} , \sigma).  \]

 \end{lemma}

 \begin{proof}

 To compute the connection one-form in mean curvature gauge, we need a normal vector $I$ that is perpendicular to $H_0$. Suppose $I$ is of the form
 \[I=I^{(0)}+\sum_{l=1}^\infty r^{-k} I^{(-k)}\] with $I^{(0)}=(1,0,0,0)$. The condition $\langle I, \partial_a \hat{X}\rangle=0$ implies
\[\langle I^{(0)}, \partial_a \hat{X}^{(0)}\rangle+\langle I^{(-1)}, \partial_a\tilde{X}\rangle=0,\] and $\langle I, H_0\rangle=0$ implies
 \[\langle I^{(0)}, H_0^{(-2)}\rangle+\langle I^{(-1)}, -2\tilde{X}\rangle=0.\] Since $\tilde{X}, \partial_1\tilde{X}, \partial_2\tilde{X}$ form an orthonormal basis in $\R^3$ and we can assume $I^{(-1)}$ is perpendicular to $I^{(0)}$, it is not hard to check that
\[I^{(-1)}=\partial_c \hat{\tau}^{(0)}\tilde{\sigma}^{cb}\partial_b\tilde{X}-\frac{1}{2}\tilde{\Delta} \hat{\tau}^{(0)}\tilde{X}\] and in general
\[I^{(-k)}=\partial_c \hat{\tau}^{(-k+1)}\tilde{\sigma}^{cb}\partial_b\tilde{X}-\frac{1}{2}\tilde{\Delta} \hat{\tau}^{(-k+1)}\tilde{X}+G( \hat{\tau}^{(0)}, \cdots, \hat{\tau}^{(-k+2)})\] for $k\geq 1$.

 The connection in mean curvature gauge is thus

\[\langle \nabla_a \frac{I}{\sqrt{-\langle I, I\rangle}}, \frac{H_0}{|H_0|}\rangle=\frac{1}{\sqrt{-\langle I, I\rangle}{|H_0|}}\langle \nabla_a I, H_0\rangle.\]

We compute
\[\begin{split}& \langle \nabla_a I, H_0 \rangle =-2r^{-2}\langle \partial_a I^{(-1)}, \tilde{X}\rangle+\sum_{k=3}^\infty r^{-k}(-2\langle \partial_a I^{(-k+1)}, \tilde{X}\rangle+\cdots+\langle \partial_a I^{(-1)}, H_0^{(-k+1)}\rangle)\\
&=r^{-2}\partial_a(\widetilde{\Delta} \hat{\tau}^{(0)}+2\hat{\tau}^{(0)})+\sum_{k=3}^\infty r^{-k}[\partial_a(\widetilde{\Delta} \hat{\tau}^{(-k+2)}+2\hat{\tau}^{(-k+2)})+G( \hat{\tau}^{(0)},\cdots, \hat{\tau}^{(-k+3)})].\end{split}\]

Since the leading term of $|H_0|$ is $\frac{2}{r}$, the leading term of $\langle I, I\rangle$ is $-1$, and the leading term of $div_{\Sigma_r}\tilde{\alpha}$ is $r^{-2}\widetilde{div}\tilde{\alpha}$ for a one-form $\tilde{\alpha}$ on $S^2$, we obtain the desired expansion for $div_{\Sigma_r} V_0$.

 \end{proof}
%CH4.2
%CH4.2
%CH4.2
\subsection{Boost $\hat{X}_r$ in $\R^{3,1}$}
Suppose $\Sigma_r$ is a family of spacelike 2-surfaces in spacetime such that

\noindent (1) The induced metric satisfies $\sigma_r = r^2 \tilde \sigma+O(r)$.

\noindent(2) The norm of the mean curvature vector satisfies $|H| =\frac{2}{r}+\frac{h^{(-2)}}{r^2} + O(r^{-3})$.

\noindent (3) The connection one-form in mean curvature gauge $V$ satisfies $ div_{\Sigma_r} V =  \frac{v^{(-3)}}{r^3} + O(r^{-4})$.

These assumptions hold on coordinate spheres of an asymptotically flat hypersurface as well as  the $r$ level surfaces at a retarded time in Bondi-Sachs coordinates. Altogether they guarantee the limit of the quasilocal energy-momentum $(e, p_1, p_2, p_3)$ with respect to isometric embeddings of $\sigma^r_{ab}$ into $\R^3$ is well-defined.

We assume the family of isometric embeddings $X_r:\Sigma\times [r_0, \infty)\rightarrow \R^{3,1}$ is of the form
\begin{equation}\label{X_r} X_r=B_r \hat{X}_r\end{equation} where $\hat{X}_r$ is a family of isometric embeddings of $\sigma_r$ into $\R^{3,1}$ that is analytic in $r$ and $B$ is a family of elements in $SO(3,1)$ that is analytic in $r$. We assume that $B=\lim_{k\rightarrow \infty} B_{k}$ where  \[B_{k}=e^{\frac{1}{r^k} b^{(-k)}}\cdots e^{\frac{1}{r} b^{(-1)}}e^{b^{(0)}}, k=0,1,2,\cdots\] for $b^{(-k)}\in \mathfrak{so}(3,1)$, the Lie algebra of $SO(3,1)$. We assume $\hat{X}_r=r\hat{X}^{(1)}+\sum_{k=0}^{\infty} r^{-k} \hat{X}^{(-k)}$
 and the time function of the embedding $\hat{X}_r$ is given by
\[\sum_{k=0}^{\infty} r^{-k} \hat{\tau}^{(-k)}.\]
Therefore we may assume $\hat{X}^{(1)}=(0, \tilde{X}_1, \tilde{X}_2, \tilde{X}_3$) is a standard embedding of $S^2$ into $\R^3$. We shall show that all $\hat{\tau}^{(0)},\cdots, \hat{\tau}^{(-k)}$ and $b^{(0)},\cdots, b^{(-k)}$ can be solved inductively. To prepare for the induction, we compute the corresponding terms in the optimal isometric embedding for $X_r$ of the given form \eqref{X_r}.

Fix an $l$, we denote
\begin{equation}\label{X_l+1} X_{l+1}=e^{r^{-l} b^{(-l)}}\cdots e^{{r}^{-1} b^{(-1)}}e^{b^{(0)}}\hat{X}.\end{equation}
Suppose $\hat{\tau}^{(0)}\cdots \hat{\tau}^{(-l)}$ and all $b^{(0)}, \cdots b^{(-l)}$ are known, we see that other components of $\hat{X}^{(0)}\cdots \hat{X}^{(-l)}$ are known and by \eqref{X_l+1} all $(X_{l+1})^{(0)}, \cdots (X_{l+1})^{(-l)}$ are known.

It is not hard to see that $X_{l+2}^{(-l)}=X_{l+3}^{(-l)}=\cdots$ is stabilized at the $r^{-l}$ term and thus $X^{(-l)}=X_{l+2}^{(-l)}$. By definition, $X_{l+2}$ is
\[\begin{split}&\{I+\cdots+r^{-l-1}[b^{(-l-1)}+G(b^{(0)},\cdots, b^{(-l)})]+\cdots\}(rB\hat{X}^{(1)}+\cdots+r^{-l} B\hat{X}^{(-l)}+\cdots)\\&=rB\hat{X}^{(1)}+\cdots +r^{-l}(b^{(-l-1)}B\hat{X}^{(1)}+B\hat{X}^{(-l)}+G(b^{(0)},\cdots, b^{(-l)}, \hat{\tau}^{(0)},\cdots, \hat{\tau}^{(-l+1)}))+O(r^{-l-1}),\end{split}\] and thus
\[X_{l+2}^{(1)}=B\hat{X}^{(1)},\]
\[X_{l+2}^{(-l)}=b^{(-l-1)}B\hat{X}^{(1)}+B\hat{X}^{(-l)}+G(b^{(0)},\cdots, b^{(-l)}, \hat{\tau}^{(0)},\cdots, \hat{\tau}^{(-l+1)}),\]
and $X_{l+1}^{(m)}$ for $0\geq m\geq -l+1$ is of the form $G( b^{(0)},\cdots, b^{(-l)}, \hat{\tau}^{(0)},\cdots, \hat{\tau}^{(-l+1)})$.

Denote by $\tau$ the time function of $X_{l+2}$. We shall plug  $\tau$ into the optimal equation, find $b^{(-l-1)}$ so that $\hat{\tau}^{(-l-1)}$ is solvable, and then solve for  $\hat{\tau}^{(-l-1)}$. Denote by $\mathfrak{G}^{(k)}$ a term of order $r^{k}$ whose coefficients depend only on $\hat{\tau}^{(0)}, \cdots, \hat{\tau}^{(-l+1)}$  and $b^{(0)}, \cdots, b^{(-l)}$ and the physical data $\sigma, |H|$, and $V$.

\begin{lemma} Write $e^{b^{(0)}}=(B_{\alpha\beta})$, then $\tau=(X_{l+2})_0$ has the following expansion:
\begin{equation}\label{tau_exp}\tau=r\tau^{(1)}+r^{-l} \tau^{(-l)}+\mathfrak{G}^{(0)}+\cdots+\mathfrak{G}^{(-l+1)}+O(r^{-l-1})\end{equation} where \[\tau^{(1)}=\sum_{i=1}^3 B_{0i}\tilde{X}_i\]  and \begin{equation}\tau^{(-l)}=\sum_{i=1}^3\sum_{\alpha=0}^3 b_{0\alpha}^{(-l-1)} B_{\alpha i}\tilde{X}_i+B_{00} \hat{\tau}^{(-l)}+G( b^{(0)},\cdots, b^{(-l)}, \hat{\tau}^{(0)},\cdots, \hat{\tau}^{(-l+1)})\end{equation}
\end{lemma}
\begin{proof}
We compute \[(B\hat{X}^{(1)})_0=B_{0\beta}(\hat{X})^{(1)}_\beta\] and recall $(\hat{X})_0^{(1)}=0$ and $(\hat{X})_i^{(1)}=\tilde{X}_i$ for $i=1, 2, 3$. Likewise, \[(b^{(-l-1)}B\hat{X}^{(1)})_0=b^{(-l-1)}_{0\alpha}B_{\alpha\beta}(\hat{X})_\beta^{(1)}.\]
 On the other hand $(B\hat{X}^{(-l)})_0=B_{00} \hat{\tau}^{(-l)}+B_{0i}\hat{X}_i^{(-l)}$, and we already know that $\hat{X}_i^{(-l)}$ depends only on $\hat{\tau}^{(-l+1)}$.
\end{proof}

Now we proceed to calculate the terms in the optimal isometric embedding equation. For simplicity, we denote $\sum_{i=1}^3\sum_{\alpha=0}^3 b_{0\alpha}^{(-l-1)} B_{\alpha i}\tilde{X}_i$ by $g$.
\begin{lemma}
With $\tau$ given by \eqref{tau_exp}, we have
\[\Delta \tau=-2r^{-1}\tau^{(1)}-2r^{-l-2}g+\mathfrak{G}^{(-2)}+\cdots +\mathfrak{G}^{(-l-2)}+O(r^{-l-3})\] and

\[\begin{split}|\nabla\tau|^2=\sum_i c_i^2-(\tau^{(1)})^2+2r^{-l-1}(\sum_i c_i d_i-\tau^{(1)}g)+\mathfrak{G}^{(-1)}+\cdots +\mathfrak{G}^{(-l-1)}+O(r^{-l-2})\end{split}\]
where
\begin{equation}\label{c_i}
c_i=B_{0i}
\end{equation}
\begin{equation}\label{d_i}
d_i=\sum_{\alpha=0}^3B_{\alpha i}b_{0 \alpha} ^{(-l-1)}
\end{equation}
\end{lemma}

\begin{proof}
We use the formula that if two functions $A$ and $B$ on $S^2$ are given by $A=A_i\tilde{X}_i$ and $B=B_j\tilde{X}_j$, then $\widetilde{\nabla}A\cdot\widetilde{\nabla}B=\sum A_i B_i-AB$.
\end{proof}

Recalling the function $f$ defined in \eqref{form_f}, we compute the expansion of $f$ and the expansions of terms that appear on the optimal isometric embedding equation \eqref{optimal3} in the following:
\begin{lemma}\label{opt_exp}
Suppose $|H_r|={2}{r}^{-1}+r^{-2} h^{(-2)}_r+\mathfrak{G}^{(-3)}+\cdots+\mathfrak{G}^{(-l-3)}+O(r^{-l-4})$ and $|H_0|={2}r^{-1}+r^{-2} h^{(-2)}_0+\mathfrak{G}^{(-3)}+\cdots+\mathfrak{G}^{(-l-3)}+O(r^{-l-4})$ then
\[f=r^{-2} f^{(-2)}+r^{-l-3}f^{(-l-3)}+\mathfrak{G}^{(-3)}+\cdots +\mathfrak{G}^{(-l-3)}+O(r^{-l-4})\] where
\begin{equation}\label{f-2} f^{(-2)}=\frac{h^{(-2)}_r-h^{(-2)}_0}{(1+|c|^2)^{1/2}}\end{equation} and
\begin{equation}\label{f-3}f^{(-l-3)}=-\frac{(h^{(-2)}_r-h^{(-2)}_0)\sum_i c_i d_i}{(1+|c|^2)^{3/2}}=f^{(-2)}(-\frac{\sum_i c_i d_i}{1+|c|^2}).\end{equation}

\[\begin{split}\Delta \sinh^{-1}\frac{(\Delta\tau)f}{|H_r||H_0|}&=\frac{1}{4} r^{-3}\widetilde{\Delta}[(\widetilde{\Delta}\tau^{(1)}) f^{(-2)}]+\frac{1}{4} r^{-l-4}\widetilde{\Delta} (f^{(-l-3)}\widetilde{\Delta} \tau^{(1)}+f^{(-2)}\widetilde{\Delta}g)\\
&+\mathfrak{G}^{(-4)}+\cdots +\mathfrak{G}^{(-l-4)}+O(r^{-l-5})\end{split}\]

\[\begin{split} div(f\nabla\tau)&=r^{-3}\widetilde{div}(f^{(-2)}\widetilde{\nabla}\tau^{(1)})+r^{-l-4}[\widetilde{div}(f^{(-2)}\widetilde{\nabla}g)
+\widetilde{div}(f^{(-l-3)}\widetilde{\nabla}\tau^{(1)})]\\
&+\mathfrak{G}^{(-4)}+\cdots +\mathfrak{G}^{(-l-4)}+O(r^{-l-5})\end{split}\]

\end{lemma}

%CH4.3
%CH4.3
%CH4.3
%CH4.3
%CH4.3
%CH4.3
\subsection{Solving the optimal isometric embedding equation for all orders}

\begin{theorem} \label{leading_order}
Suppose $\Sigma_r$ satisfies (1), (2), and (3) and the limit of the quasilocal energy-momentum $(e,p_1,p_2,p_3)$  is timelike.
There is a function $\tau^{(0)}$ on $S^2$ such that isometric embeddings $X_r:\Sigma\rightarrow \mathbb{R}^{3,1}$ with the time function $\tau$ given below solves equation (\ref{optimal}) up to $O(r^{-3})$.
\[  \tau =  (\sum_{i=1}^{3} c_{i} \tilde X_i )r + \tau^{(0)}+ O(r^{-1}) \]
where $(c_{1}, c_{2}, c_{3})$ satisfy $\frac{c_{i}}{ \sqrt{1+|c|^2}} = \frac{p_i}{e}$.
\end{theorem}

\begin{proof}
Under the assumption, the energy momentum vector $(e,p_1,p_2,p_3)$ and given by   \cite[Equation (2.6)]{wy3}
\begin{equation}\label{energy_momentum}  \int_{S^2} (h_0^{(-2)}-h^{(-2)}) dS^2 = 8 \pi e  \text{ and }  -\int_{S^2}(v_0^{(-3)}- v^{(-3)}) \tilde X_i dS^2 = 8 \pi p_i. \end{equation}

From Lemma \ref{opt_exp}, the coefficient of  the leading $r^{(-3)}$ term of $div(f\nabla\tau)-\Delta \sinh^{-1}\frac{(\Delta\tau)f}{|H_r||H_0|}$ is \[\widetilde{div}(f^{(-2)}\widetilde{\nabla}\tau^{(1)})-\frac{1}{4} \widetilde{\Delta} (f^{(-2)}\widetilde{\Delta} \tau^{(1)}),\] while the $r^{(-3)}$ term of $div V_r-div V_0$ is
\[v^{(-3)}-\frac{1}{2}\widetilde{\Delta}(\widetilde{\Delta}+2)\hat{\tau}^{(0)}.\]
Thus $\hat{\tau}^{(0)}$ is solvable if
\[\int_{S^2}[\widetilde{div}(f^{(-2)}\widetilde{\nabla}\tau^{(1)})-\frac{1}{4} \widetilde{\Delta} (f^{(-2)}\widetilde{\Delta} \tau^{(1)})-v^{(-3)}]\tilde{X}_i  dS^2\]
is zero. Computing using $\widetilde{\nabla}\tau^{(1)}\widetilde{\nabla}\tilde{X}_i=c_i-\tau^{(1)}\tilde{X}_i$, we find this is equivalent to
\[ \int_{S^2} (c_i f^{(-2)}+v^{(-3)}\tilde{X}_i) dS^2=0.\] From \eqref{energy_momentum} we see that $\int_{S^2} f^{(-2)}dS^2=\frac{-e}{\sqrt{1+|c|^2}} 8\pi$ and $\int_{S^2} v^{(-3)} \tilde{X}_i dS^2=8\pi p_i$. Thus we can take
\[  \frac{c_i}{ \sqrt{1+|c|^2}} = \frac{p_i}{e}\] and $\hat{\tau}^{(0)}$ is solvable. By equation (\ref{c_i}), $b^{(0)}$ is determined.
\end{proof}

Now we solve the optimal embedding equation to all order of $r$ by induction.

\begin{theorem}\label{all_order} Under the same assumption as in Theorem \ref{leading_order}.
Suppose the family of optimal isometric embeddings $X_r$ into $\R^{3,1}$ is of the form \eqref{X_r} described in the previous subsection. There exists $b^{(-k)}\in \mathfrak{so}(3,1)$ such that all $\hat{\tau}^{(-k)}$ can be solved for $k\geq 0$.
\end{theorem}

\begin{proof}
From Lemma \ref{opt_exp}, the coefficient of  the $r^{(-l-4)}$ term of $div(f\nabla\tau)-\Delta \sinh^{-1}\frac{(\Delta\tau)f}{|H_r||H_0|}-div(V- V_0)$ is
\[\widetilde{div}(f^{(-2)}\widetilde{\nabla}g+f^{(-l-3)}\widetilde{\nabla}\tau^{(1)})-\frac{1}{4} \widetilde{\Delta} (f^{(-l-3)}\widetilde{\Delta} \tau^{(1)}+f^{(-2)}\widetilde{\Delta}g)+\frac{1}{2}\widetilde{\Delta}(\widetilde{\Delta}+2)\hat{\tau}^{(-l-1)}+\mathfrak{G}\]
where $\mathfrak{G}$ is a term that depend on $\hat{\tau}^{(0)}, \cdots \hat{\tau}^{(-l)}$, $b^{(0)}, \cdots b^{(-l)}$, $\sigma$, $V$ and $H$.

Thus the solvability depends on whether, {for each $i=1, 2, 3$}, we can make the following expression equal to $-\int_{S^2}\mathfrak{G} \tilde X_i \, dS^2$ by choosing suitable $d_j$ in $g$:
\[\int_{S^2}\left[ \widetilde{div}(f^{(-2)}\widetilde{\nabla}g+f^{(-l-3)}\widetilde{\nabla}\tau^{(1)})-\frac{1}{4} \widetilde{\Delta} (f^{(-l-3)}\widetilde{\Delta} \tau^{(1)}+f^{(-2)}\widetilde{\Delta}g) \right]\tilde{X}_i  \, dS^2.\]

We integrate by parts and compute that this expression is equal to
\[\begin{split}&-\int_{S^2} f^{(-2)}(\tilde{X}_ig+\widetilde{\nabla}\tilde{X}_i\cdot\widetilde{\nabla}g)
+f^{(-l-3)}(\tilde{X}_i\tau^{(1)}+\widetilde{\nabla}\tilde{X}_i\cdot\widetilde{\nabla}\tau^{(1)}) \, dS^2 \\
&=-\int_{S^2} f^{(-2)}d_i+f^{(-l-3)}c_i \, dS^2\end{split}\] where we use $\widetilde{\nabla}g\widetilde{\nabla}\tilde{X}_i=d_i-g\tilde{X}_i$.

This term is \[(\int_{S^2} f^{(-2)} \,dS^2 )\sum_{j}\left(\delta_{ij}-\frac{c_i c_j}{1+\sum_k c_k^2}\right)d_j \]

Since the energy component is positive and the matrix $\delta_{ij}-\frac{c_i c_j}{1+\sum_k c_k^2}$ is positive definite, we can choose $d_j$ so that $\forall i$.
\[(\int_{S^2} f^{(-2)} \,dS^2 )\sum_{j}\left(\delta_{ij}-\frac{c_i c_j}{1+\sum_k c_k^2}\right)d_j = -\int_{S^2}\mathfrak{G} \tilde X_i \, dS^2 \]
\end{proof}

\subsection{Locally energy minimizing at all orders}

In this subsection, we show the solution $X_r=B_r\hat{X}_r$  obtained in Theorem \ref{all_order} is locally energy-minimizing at all orders. By this we mean that the second variation of the quasilocal energy is positive if we vary any $\hat{\tau}^{(-k)}$ or $b^{(-k)}$. Let $\tau$ be the time function of $X_r$ and let $\delta\tau$ be a variation. From Proposition 6.2 in \cite{wy2}, the first variation of the  energy is
\[\int_{{\Sigma_r}}\left[-(\widehat{H} \hat{\sigma}^{ab} -\hat{\sigma}^{ac} \hat{\sigma}^{bd} \hat{h}_{cd})\frac{\nabla_b\nabla_a \tau}{\sqrt{1+|\nabla\tau|^2}}+ div_{\Sigma_r} (\frac{\nabla\tau}{\sqrt{1+|\nabla\tau|^2}} \cosh\theta|{H}|-\nabla\theta-V_r) \right] \delta\tau d\Sigma_r. \]

By the derivation in \S\ref{opt}, we can rewrite this as
\[ \int_{\Sigma_r} \left[ div(f\nabla\tau)-\Delta \sinh^{-1}\frac{(\Delta\tau)f}{|H_r||H_0|}-div(V- V_0) \right] \delta \tau d\Sigma_r.\]

From Lemma \ref{opt_exp}, the coefficient of  the $r^{-l-4}$ term of $div(f\nabla\tau)-\Delta \sinh^{-1}\frac{(\Delta\tau)f}{|H_r||H_0|}-div(V- V_0)$ is
\[\begin{split}&\widetilde{div}(f^{(-2)}\widetilde{\nabla}g+f^{(-l-3)}\widetilde{\nabla}\tau^{(1)})-\frac{1}{4} \widetilde{\Delta} (f^{(-l-3)}\widetilde{\Delta} \tau^{(1)}+f^{(-2)}\widetilde{\Delta}g)\\
&+\frac{1}{2}\widetilde{\Delta}(\widetilde{\Delta}+2)\hat{\tau}^{(-l-1)}+{G}( b^{(0)}, \cdots b^{(-l)}, \hat{\tau}^{(0)}, \cdots \hat{\tau}^{(-l)})\end{split}\] where $g=\sum_{i=1}^3\sum_{\alpha=0}^3 b_{0\alpha}^{(-l-1)} B_{\alpha i}\tilde{X}_i$.

On the other hand, the $r^{-l}$ term of $\tau$ is
\begin{equation}\tau^{(-l)}=g+B_{00} \hat{\tau}^{(-l)}+G( b^{(0)},\cdots, b^{(-l)}, \hat{\tau}^{(0)},\cdots, \hat{\tau}^{(-l+1)}).\end{equation}

When we consider the variation of  $\hat{\tau}^{(-l)}$, $\delta\hat{\tau}^{(-l)}$, the leading term term of the second variation is of the order $r^{-2l-3}$ with coefficient
\[\frac{B_{00}}{2}\int_{S^2} [\widetilde{\Delta}(\tilde{\Delta}+2)\delta \hat{\tau}^{(-l)}]\delta \hat{\tau}^{(-l)} dS^2.\]
We may assume $\int_{S^2} \delta \hat{\tau}^{(-l)} dS^2=0$ by normalization.
By decomposing $\delta \hat{\tau}^{(-l)}$ into sum of eigenfunctions of $S^2$ and noting that the first non-zero eigenvalue of $S^2$ is $-2$, this is always positive.

Varying $b^{(-l-1)}$ is equivalent to varying $g$ and $\delta g=\sum_{i=1}^3\sum_{\alpha=0}^3 \delta b_{0\alpha}^{(-l-1)} B_{\alpha i}\tilde{X}_i$.
The coefficient of the leading $r^{-2l-4}$ term of the second variation with respect to $b^{(-l-1)}$ is then
\[\int_{S^2} [\widetilde{div}(f^{(-2)} \widetilde\nabla \delta g) - \frac{\widetilde \Delta(f_2\widetilde \Delta \delta g)}{4}    ]\delta g dS^2.\]
Using integration by part, we derive that this is equal to
\[-\int_{S^2}  f^{(-2)} [ |\widetilde\nabla \delta g|^2 + (\delta g )^2]  dS^2.\]
This simply gives a positive multiple of quasilocal energy since
\[ |\widetilde\nabla \delta g|^2 + (\delta g )^2 = \sum_{i}(\sum_\alpha   \delta b_{0\alpha}^{(-l-1)} B_{\alpha i})^2 \] is a positive constant and $-\int_{S^2} f^{(-2)} dS^2=\frac{8\pi e}{\sqrt{1+|c|^2}}$.
 
\end{document}